\documentclass{ifacconf}

\usepackage{natbib}        
\usepackage{amsfonts}
\usepackage{cite}
\usepackage{amssymb}
\usepackage{xcolor}
\usepackage{mathtools}
\usepackage{dsfont}
\usepackage{mathrsfs}
\usepackage{epsfig}
\usepackage{caption}
\usepackage{nicematrix}
\usepackage{algorithm}
\usepackage{algpseudocode}
\usepackage{tabularx, booktabs, array}
\newtheorem{ass}{Assumption}

\newcommand{\argmin}{\mathop{\rm argmin}}

\newcommand{\trace}[1]{\text{Tr}\left(#1\right)}

\newcommand{\Ex}[2][]{\mathds{E}_{#1} \left[ #2 \right]}

\newcommand{\Sset}{\mathbb{S}}

\newcommand{\Bcal}{{\cal B}}

\newcommand{\Hcal}{{\cal H}}
\newcommand{\Kcal}{{\mathcal{K}}}

\newcommand{\Mcal}{{\cal M}}
\newcommand{\Ncal}{{\mathcal{N}}}
\newcommand{\Pcal}{{\cal P}}

\newcommand{\Rcal}{{\cal R}}

\newcommand{\Tcal}{{\mathcal{T}}}

\newcounter{l1}
\newcounter{l2}
\newcounter{l3}
\newcommand{\bdotlist}{\begin{list}{$\bullet$}{}}
\newcommand{\bboxlist}{\begin{list}{$\Box$}{}}
\newcommand{\bbboxlist}{\begin{list}{\raisebox{.005in}{{\tiny
$\blacksquare$ \ \ }}}{}}
\newcommand{\bdashlist}{\begin{list}{$-$}{} }
\newcommand{\blist}{\begin{list}{}{} }
\newcommand{\barablist}{\begin{list}{\arabic{l1}}{\usecounter{l1}}}
\newcommand{\balphlist}{\begin{list}{(\alph{l2})}{\usecounter{l2}}}
\newcommand{\bAlphlist}{\begin{list}{\Alph{l2}.}{\usecounter{l2}}}
\newcommand{\bdiamlist}{\begin{list}{$\diamond$}{}}
\newcommand{\bromalist}{\begin{list}{(\roman{l3})}{\usecounter{l3}}}

\newcommand{\dist}{\mathbb{P}}

\begin{document}
\begin{frontmatter}

\title{Distributional Robustness in Output Feedback Regret-Optimal Control\thanksref{footnoteinfo}} 

\thanks[footnoteinfo]{Funded by Deutsche Forschungsgemeinschaft (DFG, German Research Foundation) under Germany's Excellence Strategy - EXC 2075 – 390740016. We acknowledge the support by the Stuttgart Center for Simulation Science (SimTech).}

\author[First]{Shuhao Yan} 
\author[First]{Carsten W. Scherer}

\address[First]{Department of Mathematics, University of Stuttgart, Germany (e-mail: \{shuhao.yan, carsten.scherer\}@imng.uni-stuttgart.de).}

\begin{abstract}                
This paper studies distributionally robust regret-optimal (DRRO) control with purified output feedback for linear systems subject to additive disturbances and measurement noise. These uncertainties (including the initial system state) are assumed to be stochastic and distributed according to an unknown joint probability distribution within a Wasserstein ambiguity set. We design affine controllers to minimise the worst-case expected regret over all distributions in this set. The expected regret is defined as the difference between an expected cost incurred by an affine causal controller and the expected cost incurred by the optimal noncausal controller with perfect
knowledge of the disturbance trajectory at the outset. 
Leveraging the duality theory in distributionally robust optimisation, we derive strong duality results for worst-case expectation problems involving general quadratic objective functions, enabling exact reformulations of the DRRO control problem as semidefinite programs (SDPs). Focusing on one such reformulation, we eliminate certain decision variables. This technique also permits a further equivalent reformulation of the SDP as a distributed optimisation problem, with potential to enhance scalability.
\end{abstract}


\end{frontmatter}

\section{INTRODUCTION}
Distributionally robust optimisation (DRO) is a powerful model for optimisation under uncertainty, incorporating robust and stochastic optimisation as special cases. It seeks decisions to minimise an expected cost for all probability distributions in a pre-specified set, referred to as the \emph{ambiguity set}. In particular, ambiguity sets defined in terms of the Wasserstein metric have been shown to provide desirable out-of-sample performance guarantees \citep{mohajerin2018data}. In DRO, the inner maximisation of an expected cost over all distributions in an ambiguity set is generally infinite-dimensional and intractable, motivating considerable efforts towards efficient optimisation and tractable reformulations of the worst-case expectation, see, e.g. \citep{Kuhn19DROtutorial}.

In addition to using conventional metrics, such as cost, there has been increasing interest in considering regret as a desirable alternative, for example, in estimator design \citep{eldar2004competitive} and adaptive control \citep{dean2018regret}. More recently, \citet{9483023} and \citet{10061197} develop regret-optimal control for linear quadratic problems, where regret-optimal control problems can be reduced to $H_\infty$ control problems. Inspired by their work, \citet{10384311} propose distributionally robust regret-optimal (DRRO) control, and \citet{martin2024regret} extend the regret-optimal control framework to consider linear time-varying uncertain systems subject to both parametric uncertainty and additive disturbances.

In our previous work \citep{10384311}, we design linear controllers to minimise the worst-case expected regret over all distributions in a Wasserstein ball. As demonstrated by an example, such controllers can outperform distributionally robust controllers, which minimise the worst-case expected cost.
In this paper, we extend the DRRO control framework to consider output feedback settings, where true outputs are corrupted by additive noise, and parameterise the input trajectory $u$ as an affine function of the purified output trajectory $\eta$, i.e. $u=K\eta + g$ where $K$ is a feedback gain matrix and $g$ is a vector. \citet{10313386} also investigate DRRO control with output feedback, while they use a linear control parameterisation and an alternative benchmark policy in the definition of the expected regret. After reformulating the control problem using \citep[Theorem 2]{10384311}, the control synthesis requires solving a sequence of semidefinite programs (SDPs) due to optimising over variables alternately.

\textit{Summary of contributions:} In this paper, we consider two types of Wasserstein ambiguity sets, denoted by $\Bcal_1$ and $\Bcal_2$, to capture possible distributional uncertainty affecting the system. The set $\Bcal_1$ is centred at a probability distribution that is absolutely continuous with respect to the Lebesgue measure, and $\Bcal_2$ is centred at an empirical distribution. Based on each set of assumptions regarding the Wasserstein ball, we derive a corresponding strong duality result for worst-case expectation problems with general quadratic objectives. Our results complement several related existing results, i.e. \citep[Theorem 11]{Kuhn19DROtutorial}, \citep[Theorem 2]{10384311}, and \citep[Lemma 3]{brouillon2023distributionally}, on exact reformulations of such worst-case expectation problems. We leverage the derived duality results to reformulate the DRRO control problem equivalently as SDPs and then focus on the SDP reformulation where the Wasserstein ball $\Bcal_1$ is assumed. In this SDP, we eliminate the affine term, $g$, of the chosen input parameterisation to reduce computational complexity and show it admits an explicit expression. Exploiting the Projection Lemma \citep[Lemma 3.1]{gahinet1994linear} and its generalisation \citep[Theorem 2]{SCHERER1995327} to block triangular matrix variables, we can further reduce the number of decision variables of this SDP by eliminating the feedback gain matrix $K$ of the input parameterisation at the cost of additional constraints. We show this technique can reduce computation time in Section \ref{sec:numerical example}. More importantly, eliminating these variables yields an SDP that can be rewritten equivalently as distributed optimisation with consensus conditions, enabling the application of tailored distributed algorithms and thus the handling of large-scale problems. Such a distributed reformulation is not possible with the original SDP.

\textit{Notation:} Given vectors $v_1, \ldots, v_n$, $\mbox{col}(v_1, \ldots, v_n)$ denotes their vertical concatenation; $\|v_1\|$ denotes the Euclidean norm of $v_1$; and $\|v_1\|^2_P:=v_1^\top P v_1$, where $P$ is a square matrix. 
For a symmetric block matrix, we use $*$ to represent some sub-block that can be inferred from symmetry. Let $\Mcal(\Xi)$ denote the set of Borel probability measures supported on $\Xi$ with finite second moments. 

\section{PROBLEM FORMULATION}

\subsection{System Description}
Consider an uncertain linear time-invariant (LTI) system evolving over a finite horizon of length $T$, described by
\begin{equation} \label{eq:LTI}
    x_{t+1} = A x_t + B u_t + w_t, ~~ y_t =H x_t + v_t, 
\end{equation}
where $x_t \in \Rset^{n_x}$, $u_t \in \Rset^{n_u}$ and $y_t \in \Rset^{n_y}$. We assume the initial system state $x_0$, the disturbances $w_0,\ldots,w_{T-1}$ and the measurement noise $v_0,\ldots, v_{T-1}$ are random variables distributed according to a joint distribution $\dist^\star$. This distribution is unknown but assumed to lie in a given compact set $\Pcal$. We define finite horizon trajectories as
$x:=\mbox{col}(x_0, \dots, x_T) \in \Rset^{N_x}$, $u\!:=\! \mbox{col}(u_0, \dots, u_{T-1}) \in \Rset^{N_u}$, $w\!:=\! \mbox{col}(x_0, w_0, \dots, w_{T-1}) \!\in \Rset^{N_x}$, $y:=\mbox{col}(y_0,\ldots,y_{T-1}) \in \Rset^{N_y}$, and $v:=\mbox{col}(v_0,\ldots,v_{T-1}) \in \Rset^{N_y}$, where $N_x:=(T+1)n_x$, $N_u:=Tn_u$, and $N_y:=Tn_y$. 
On a finite horizon, the system \eqref{eq:LTI} can then be expressed in terms of linear causal mappings between trajectories as
\begin{equation}\label{eq:compact dynamics}
        x = F u + G w, ~~ y = C x +v, 
\end{equation}
where $F$ and $G$ are block lower triangular Toeplitz matrices and $C$ is a rectangular diagonal matrix. They can be routinely constructed from matrices $A, B$, and $H$.

\subsection{Controller Parameterisation}
We parameterise control inputs as causal affine mappings from \emph{purified} outputs \citep{ben2006extending}. To this end, we first introduce a fictitious LTI system without disturbances and measurement noise, described by
\begin{equation}\label{eq:fictitious system}
       \hat{x}_0=0, ~~ \hat{x}_{t+1} =A \hat{x}_t + Bu_t, ~~\hat{y}_t =H \hat{x}_t. 
\end{equation}
A purified output is defined as $\eta_t:=y_t -\hat{y}_t$. Finite horizon trajectories corresponding to \eqref{eq:fictitious system} are defined as 
$\hat{x}:=\mbox{col}(\hat{x}_0,\ldots,\hat{x}_T)\in\Rset^{N_x}$, $\hat{y}:=\mbox{col}(\hat{y}_0, \ldots,\hat{y}_{T-1}) \in\Rset^{N_y}$, and $\eta:=\mbox{col}(\eta_0,\ldots,\eta_{T-1})\in\Rset^{N_y}$.
Then, the controller parameterisation is given by 
\begin{equation}
    u=K\eta + g,   \label{eq:purified output feedback}
\end{equation}
where $g\in\Rset^{N_u}$ and
\vspace{-3.5mm}
\[
K\!:=\!
\begin{pmatrix}
    K_{0,0} &\!  0 & \cdots & 0 \\
    K_{1,0} &\! K_{1,1} & \cdots &  0\\
    \vdots &\! \vdots & \ddots & \vdots\\
    K_{T\!-\!1,0} &\! \cdots & \cdots & K_{T\!-\!1,T\!-\!1}
\end{pmatrix}\!\!=\!\sum_{i=1}^T L_i \!\overbrace{\begin{pmatrix}
    K_{i-1,i-1}\\
    K_{i,i-1}\\
    \vdots\\
    K_{T\!-\!1,i-1}
\end{pmatrix}}^{K_i:=}\! R_i.
\]
The left and right factors $L_i$ and $R_i$ are straightforward to construct and not displayed explicitly for reasons of space. The block matrix between them is denoted by $K_i$.
We use $\Kcal$ to denote the subset of $\Rset^{N_u\times N_y}$ that consists of all block lower triangular matrices that correspond to this parameterisation\footnote{As shown in the literature, e.g. \citep{5742976,NEURIPS2023_3b7a66b2}, \eqref{eq:purified output feedback} is equivalent to the parameterisation $u=My+s$, where $M\!\in\!\Kcal$ and $s\!\in\!\Rset^{N_u}$. Still, the former may be preferred in terms of implementation, as the latter requires a parameter transformation to ensure linear dependence of $x$ and $u$ on decision variables and then possibly large-scale matrix inversion
to recover $(M,s)$.}. Based on \eqref{eq:purified output feedback}, we can express the input and state trajectories of the LTI system \eqref{eq:compact dynamics} as
\[
     u= K(CGw+v) + g, ~ x= FK (CGw +v) +Fg + Gw,
\] 
with an affine dependence on $(K,g)$. This ensures that the minimisation of a convex function of $(x,u)$ is convex.

\subsection{Control Problem}

We aim to synthesise output feedback controllers in the form of \eqref{eq:purified output feedback} that minimise the worst-case expected regret incurred by them over all distributions in the ambiguity set $\Pcal$. Considering a quadratic cost given by
\begin{equation}
    J(u,w):=x^\top J_x x + u^\top J_u u \label{eq:cost function}
\end{equation}
where $J_x\in\Sset_+^{N_x}$ and $J_u\in\Sset_{++}^{N_u}$, we can define the regret incurred by an input sequence $u$ with respect to a benchmark policy $u^\star(\cdot)$ as
\begin{equation}
    J(u,w) - J(u^\star(\cdot),w), \nonumber 
\end{equation}
where $u^\star(\cdot)$ possibly depends on $w$ or $v$.
In general, a benchmark is chosen as one that yields superior performance but is not implementable in practice. In this work, we choose the benchmark as
\begin{equation}
    u^\star(w):=\argmin_{u\in\Rset^{N_u}} J(u,w) =K^\star w, \label{eq:noncausal benchmark}
\end{equation}
where $K^\star:=-(J_u+F^\top J_x F)^{-1} F^\top J_x G$, and it has knowledge of the disturbance realisations at the outset. Other benchmark policies considered in the literature include, e.g. noncausal output feedback control that is optimal in an $\Hcal_2$ sense \citep{10313386}, and a benchmark that has knowledge of the underlying distributions \citep{10534845} and is used to define the \emph{ex-ante} regret \citep{poursoltani2024risk}. The optimality of $u^\star(w)$ in \eqref{eq:noncausal benchmark} allows us to express the regret in a quadratic form as
\begin{equation}
    J(u,w) - J(u^\star(w),w)=(u-K^\star w)^\top D (u-K^\star w) \label{eq:quadratic form of regret}
\end{equation}
for any input sequence $u\in\Rset^{N_u}$, where $D:=J_u+F^\top J_x F$. 
Therefore, based on the parameterisation \eqref{eq:purified output feedback}, we denote the regret as
$\Rcal(K,g,w,v) := J(K\eta +g,w) - J(K^\star w,w)$
\[ ~~   
    =\bigl\| (KCG\!-\!K^\star \, ~  K )
    \begin{pmatrix}
        w\\
        v
    \end{pmatrix} +g \bigr\|^2_D. 
\]
We now formulate the \emph{distributionally robust regret-optimal control} problem as
\begin{equation}\label{eq:drro control}
    \inf_{K\in\Kcal,\, g \in\Rset^{N_u}}\sup_{\dist\in\Pcal} \Ex[\dist]{\Rcal(K,g,w,v)}.
\end{equation}
To ensure this is meaningful, we assume $K^\star \neq 0$ hereafter.

\section{WORST-CASE EXPECTATION PROBLEM}\label{sec: worst-case expectation problems}
This section derives two duality results for worst-case expectation problems with general quadratic objectives by leveraging a well-established duality theory in DRO. These results facilitate equivalent reformulations of problem \eqref{eq:drro control}. 

We first introduce the 2-Wasserstein distance between probability distributions. 
\begin{defn} \label{defn:wass_dist} 
The \emph{2-Wasserstein distance} between two distributions $\dist_1, \dist_2 \in \Mcal(\Xi)$  is defined as 
\begin{align*}
    W_2(\dist_1, \, \dist_2)^2 := \hspace{-.075in}\inf_{\pi \in \Pi(\dist_1, \, \dist_2)}  \int_{\Xi \times \Xi} \hspace{-.05in} \| z_1 - z_2 \|^2 \pi(d z_1, \, d z_2),
\end{align*}
 where $\Pi(\dist_1, \, \dist_2)$ denotes the set of all joint distributions in $\Mcal(\Xi \times \Xi)$ with marginal distributions $\dist_1$ and $\dist_2$.
\end{defn}

Given a \emph{nominal distribution} $\dist_0 \in \Mcal( \Xi )$, we define the ambiguity set $\Pcal$ based on the 2-Wasserstein distance as
\begin{align} \label{eq:amb_set}
    \Pcal := \{ \dist \in \Mcal( \Xi ) \, | \, W_2(\dist, \dist_0) \leq r \}.
\end{align}
The set $\Pcal$ is also referred to as a Wasserstein ball. Its radius $r$ represents the belief of a decision maker about the size of distributional uncertainty. In the following, we make two sets of assumptions regarding the set $\Xi$ and the nominal distribution $\dist_0$, respectively.

\subsection{Unbounded $\Xi$ and Absolutely Continuous $\dist_0$ with Respect to the Lebesgue Measure}
\begin{ass}\label{assumption: support on Rn}
    The set $\Xi$ is equal to $\Rset^n$.
\end{ass}
\begin{ass}\label{ass: absolute continuity}
    The nominal distribution $\dist_0 \in \Mcal(\Rset^n)$ is absolutely continuous with respect to the Lebesgue measure on $\Rset^{n}$.
\end{ass}

\begin{thm}\label{theorem: strong duality for worst-case expectation problem with quadratic objective: absolutely continuous P0}
    Let $r>0$ and Assumptions \ref{assumption: support on Rn}-\ref{ass: absolute continuity} hold. Consider a family of worst-case expectation problems given by
    \begin{equation}
        \sup_{\dist\in\Pcal}\Ex[\dist]{\xi^\top P \xi + 2q^\top \xi + c}, \label{eq:worst-case expectation general quadratic: restriction on P}
    \end{equation}
    where $P \in \Sset^{n}$, $q\in\Rset^{n}$, $c\in\Rset$ and $P$ satisfies $\lambda_{\rm max}(P)\neq 0$. Then, the optimal value of problem \eqref{eq:worst-case expectation general quadratic: restriction on P} is finite and equal to the optimal value of the SDP:
    \begin{align*}
        &\inf_{\gamma\geq0,\, \beta\in\Rset,\, X\in\Sset^n} \gamma \left(r^2 - \trace{\Sigma_0} - \|\mu_0\|^2 \right) + \trace{X} + \beta \\
        &\begin{aligned}
            \text{ s.t. } 
            &\begin{pmatrix}
                X & \gamma \Lambda^\top\\
                \gamma \Lambda & \gamma I\!-\!P
            \end{pmatrix}\!\succeq\! 0,
            ~
            \begin{pmatrix}
                \beta-c &* \\ 
                \gamma \mu_0 + q & \gamma I-P
            \end{pmatrix}\!\succeq\! 0,
            ~\gamma I  \!\succ\! P,
        \end{aligned}     
    \end{align*}
    where $\Sigma_0$ and $\mu_0$ are the covariance and mean of nominal distribution $\dist_0$, respectively, and $\Lambda$ is the lower triangular matrix such that $\Sigma_0=\Lambda \Lambda^\top$.
\end{thm} 
Theorem \ref{theorem: strong duality for worst-case expectation problem with quadratic objective: absolutely continuous P0} generalises \citep[Theorem 2]{10384311} and follows using similar reasoning, while the latter considers a class of worst-case expectation problems with purely quadratic objectives.
Under Assumption \ref{ass: absolute continuity}, we further establish that $\Sigma_0\succ 0$  
and thus $\Lambda$ is unique.

\subsection{Ellipsoidal $\Xi$ and Finitely Supported $\dist_0$}
\begin{ass}\label{ass:ellipsoidal set with one interior point}
The set $\Xi\!:=\!\{\xi \!\in\! \Rset^n |\, \xi^\top\! P_2 \xi + 2 q_2^\top \xi + c_2 \!\leq\! 0\}$ has an interior point, where $P_2\succ 0$, $q_2\in\Rset^n$ and $c_2\in\Rset$.
\end{ass}

\begin{ass}\label{assumption:discrete nominal distribution}
    The nominal distribution $\dist_0\in\Mcal(\Xi)$ is a discrete uniform distribution with support $\{\zeta_1,\ldots,\zeta_N\}$. 
\end{ass}

\begin{thm}\label{thm:strong duality for worst-case expectation problem over ambiguity set centered on discrete distribution}
Let $r>0$ and Assumptions \ref{ass:ellipsoidal set with one interior point}-\ref{assumption:discrete nominal distribution} hold. Consider a quadratic function $\ell(\xi)=\xi^\top P_1 \xi + 2q_1^\top \xi + c_1$, where $P_1\in\Sset^n$, $q_1\in\Rset^{n}$ and $c_1\in\Rset$, the optimal value of the worst-case expectation problem 
{\setlength{\abovedisplayskip}{2pt}
 \setlength{\belowdisplayskip}{2pt}
\begin{equation}
    v_{opt}:=\sup_{\dist\in\Pcal}\Ex[\dist]{\ell(\xi)} \label{eq:quadratic worst-case expectation: ellipsoidal support}
\end{equation}}
is finite and equal to the optimal value of the SDP:
\begin{align}
    &\inf_{\lambda\geq0,\alpha_i\geq0,\gamma_i\in\Rset} \lambda r^2  \!+\! \frac{1}{N}\sum_{i=1}^N \gamma_i \label{eq:sdp with n lmis}\\
    &
    \text{\,s.t. } 
                 \Omega(\lambda,\zeta_i,\alpha_i,-\gamma_i)
                 \succeq 0,~  \forall i=1,\ldots,N, \nonumber
\end{align}
where $\Omega$ is a matrix-valued function and defined as
\[
\Omega(\lambda,\zeta,\alpha,t):=\begin{pmatrix}
        \lambda\|\zeta\|^2-c_1 +\alpha c_2 - t &*\\ 
        \alpha q_2 - \lambda \zeta -q_1 & \lambda I -P_1 + \alpha P_2
    \end{pmatrix}.
\]
\end{thm}
\vspace{-10pt}
\begin{pf}
It is clear that $(\Xi,\|\cdot\|)$ is a Polish space and the quadratic function $\ell(\cdot)$ is Borel measurable, satisfying $\Ex[\dist_0]{|\ell(\zeta)|}<\infty$. Therefore, from \citep[Theorem 1]{gao2023distributionally}, it follows that
\begin{equation}\label{eq:initial form of strong duality}
    \sup_{\dist\in\Pcal}\Ex[\dist]{\ell(\xi)}
    =
    \inf_{\lambda\geq0}\{\lambda r^2 - \Ex[\dist_0]{\phi(\lambda,\zeta)}\}, 
\end{equation}
where $\phi(\lambda,\zeta):=\inf_{\xi\in\Xi}\{\lambda \|\xi-\zeta\|^2 - \ell(\xi) \}=\inf_{\xi\in\Xi}\{\xi^\top (\lambda I -P_1)\xi -2(\lambda \zeta + q_1)^\top \xi - c_1\} + \lambda \zeta^\top \zeta$. As $\Xi$ is compact by Assumption \ref{ass:ellipsoidal set with one interior point} and $\xi^\top (\lambda I \!-\! P_1)\xi \!-\! 2(\lambda \zeta \!+\! q_1)^\top \xi \!-\! c_1$ is a continuous function of $\xi$ for all $\zeta \!\in\! \Xi$ and all $\lambda \!\in\! [0,\infty)$, 
$
    \inf_{\xi\in\Xi}\{\xi^\top (\lambda I \!-\!P_1)\xi \!-\! 2(\lambda \zeta \!+\! q_1)^\top \xi \!-\! c_1\} 
$
is finite for all $\zeta\!\in\!\Xi$ and all $\lambda\!\in\![0,\infty)$. This, along with finiteness of the second moment of $\dist_0$, shows that $\Ex[\dist_0]{\phi(\lambda,\zeta)}\!>\!-\infty$ for all $\lambda\!\in\![0,\infty)$. Thus, $v_{opt}$ is finite. 

We next show that $\phi(\lambda,\zeta)$ is equal to the optimal value of an SDP. 
Since $\phi(\lambda,\zeta)$ is the optimal value of a quadratic program with a single quadratic constraint, which is strictly feasible under Assumption \ref{ass:ellipsoidal set with one interior point}, strong duality holds for such a problem. Therefore, it holds that
\begin{equation} 
    \phi(\lambda,\zeta)=\sup_{\alpha\geq 0,\, t\in\Rset}\{t \,|\, \Omega(\lambda,\zeta,\alpha,t)\succeq 0 \} . \nonumber
\end{equation}
As $\dist_0$ is a discrete distribution, we then have that
{\setlength{\belowdisplayskip}{2pt}
\begin{align}
    &\inf_{\lambda\geq0}\{ \lambda r^2 \!-\!\Ex[\dist_0]{\phi(\lambda,\zeta)}\}
    =
    \inf_{\lambda\geq0}\biggl\{ \lambda r^2  \!-\!\frac{1}{N}\sum_{i=1}^N \phi(\lambda,\zeta_i) \biggr\} \nonumber \\[-4pt]
    =&\inf_{\lambda\geq0}\biggl\{ \lambda r^2 
    \!-\frac{1}{N} \sum_{i=1}^N \sup_{\alpha_i\geq0,t_i\in\Rset}\biggl\{t_i \left|\, 
    \Omega(\lambda,\zeta_i,\alpha_i,t_i)
    \succeq 0 
    \biggr\} \right.
    \biggr\} \nonumber \\[-3pt]
    =&
    \inf_{\lambda\geq0}\biggl\{ \lambda r^2 
    \!+\!
    \frac{1}{N} \sum_{i=1}^N \inf_{\alpha_i\geq0,\gamma_i\in\Rset}\biggl\{\gamma_i \left|\,
    \Omega(\lambda,\zeta_i,\alpha_i,-\gamma_i)
    \succeq 0
    \biggr\}\right.
    \biggr\} \nonumber \\[-3pt]
    =& \text{~\eqref{eq:sdp with n lmis}} \nonumber, 
\end{align}}%
where the penultimate equality follows from substituting $-\gamma_i$ for $t_i$. This, along with \eqref{eq:initial form of strong duality}, proves the desired result.
\end{pf}
In case $\dist_0$ is a discrete distribution with non-uniform weights $p_i$, 
Theorem \ref{thm:strong duality for worst-case expectation problem over ambiguity set centered on discrete distribution} still holds true with $p_i$ replacing ${1}/{N}$ in \eqref{eq:sdp with n lmis} and appearing inside the summation over $i$.

Theorem \ref{thm:strong duality for worst-case expectation problem over ambiguity set centered on discrete distribution} complements Theorem \ref{theorem: strong duality for worst-case expectation problem with quadratic objective: absolutely continuous P0} without imposing a condition on the quadratic coefficient matrix, as well as several related results, e.g. \citep[Lemma 3]{brouillon2023distributionally} and \citep[Theorem 11]{Kuhn19DROtutorial}. Both provide exact SDP reformulations of worst-case expectation problems over Wasserstein balls centred at empirical distributions, with the former assuming distributions are supported on a compact polytope and the latter on $\Rset^n$.

\begin{rem}
    The computational complexity of SDP \eqref{eq:sdp with n lmis} grows with $N$, while it is typically desired to use a large number of data points (if available) to obtain a smaller Wasserstein ball or to ensure that the Wasserstein ball contains the true distribution with a higher confidence.
To remedy this issue, one can 
compress an $N$-point empirical distribution to obtain an $m$-point distribution, such that the latter is close to the former in the Wasserstein distance and $m \ll N$ \citep{rujeerapaiboon2022scenario}. Then, the resulting distribution can be chosen as the central distribution to construct a suitable Wasserstein ball. 
\end{rem}

\section{DISTRIBUTIONALLY ROBUST REGRET-OPTIMAL Control}
The duality results in Section \ref{sec: worst-case expectation problems} allow us to address the inner maximisation of problem \eqref{eq:drro control}, despite its generally infinite-dimensional nature. Building on them, we reformulate \eqref{eq:drro control} as a tractable, finite-dimensional, convex problem in Section \ref{sec:convex reformulation}, followed by variable reduction of the reformulated problem in Sections \ref{sec: optimality of linear policy} and \ref{sec: complexity reduction}.

\subsection{SDP Reformulation}\label{sec:convex reformulation}
As a direct application of Theorem \ref{theorem: strong duality for worst-case expectation problem with quadratic objective: absolutely continuous P0} and the Schur complement, we present an exact reformulation of the DRRO control problem \eqref{eq:drro control} as an SDP.
\begin{thm}\label{theorem:initial sdp reformulation for drro control}
Let $r>0$ and Assumptions \ref{assumption: support on Rn}-\ref{ass: absolute continuity} hold with $n=N_x+N_y$. Then, problem \eqref{eq:drro control} is equivalent to
\begin{subequations}\label{eq:exact sdp reformulation of drro control: original}
   \begin{align}
        &\hspace{-60pt}\underset{X\in\Sset^{N_x\!+\!N_y}}{\underset{\gamma\geq0,\, \beta\in\Rset}{\underset{K\in\Kcal,\, g\in\Rset^{N_u}}{\inf}}} \gamma \left(r^2 - \trace{\Sigma_0} - \|\mu_0\|^2 \right) + \trace{X} + \beta \\[-10pt]
            \text{s.t. } &\hspace{0pt} \begin{pmatrix}
                \gamma I & *\\ 
                (KCG\!-\!K^\star \, ~ K) & D^{-1}
            \end{pmatrix} \succ 0, \label{eq: LMI 1}\\
           Q(\gamma, X, K):= &\begin{pmatrix}
                X & * & * \\ 
                \gamma \Lambda & \gamma I & *\\ 
                0 & (KCG\!-\!K^\star \, ~ K) & D^{-1}
            \end{pmatrix}\succeq 0,  \label{eq: LMI 2} \\
            &\begin{pmatrix}
                \beta & * & * \\ 
                \gamma \mu_0 & \gamma I &* \\ 
                -g &\, (KCG\!-\!K^\star \, ~ K) \!& D^{-1}
            \end{pmatrix}\succeq 0,  \label{eq: LMI 3}
    \end{align} 
    \end{subequations}
where $Q$ is a matrix-valued function.
\end{thm}
Recalling the expression for the regret $\Rcal(K,g,w,v)$, we note that $(KCG\!\!-\!\!K^\star \, ~ K)^\top D (KCG\!\!-\!\!K^\star \, ~ K)\!\succeq\! 0$ for all $K\!\in\!\Kcal$ and its maximum eigenvalue is positive as $D\!\succ\! 0$ and $K^\star \!\neq\! 0$ by assumption. 
This makes Theorem \ref{theorem: strong duality for worst-case expectation problem with quadratic objective: absolutely continuous P0} applicable to deriving such a reformulation in \eqref{eq:exact sdp reformulation of drro control: original}.    

Alternatively, under Assumptions \ref{ass:ellipsoidal set with one interior point}-\ref{assumption:discrete nominal distribution}, one can derive another exact reformulation of problem \eqref{eq:drro control} based on Theorem \ref{thm:strong duality for worst-case expectation problem over ambiguity set centered on discrete distribution}, and the corresponding reformulation constitutes a data-driven approach to controller design. We omit that reformulation to streamline the presentation of this paper. 

The complexity of SDP \eqref{eq:exact sdp reformulation of drro control: original} grows polynomially with state, input, and output dimensions and the time horizon. It becomes challenging to solve in high dimensions. To mitigate this issue, one possible option is to break the problem into two stages, where the second stage does not require optimisation. As will be introduced in the following sections, we eliminate certain decision variables from \eqref{eq:exact sdp reformulation of drro control: original} step by step to achieve an SDP with much fewer variables, and those eliminated variables can be reconstructed efficiently using standard linear algebra. 

\subsection{Elimination of the Decision Variables $g$ and $\beta$}\label{sec: optimality of linear policy}
We first eliminate variables $g$ and $\beta$ from problem \eqref{eq:exact sdp reformulation of drro control: original}.
\begin{lem}\label{lemma: eliminate g and beta}
    The optimal value of problem \eqref{eq:exact sdp reformulation of drro control: original} is equal to 
    \begin{equation}\label{problem: exact sdp reformulation of drro control: without g and beta}
        \underset{X\in\Sset^{N_x\!+\!N_y}}{\underset{{K\in\Kcal,\, \gamma\geq0}}{\inf}}\! \bigl\{ \gamma \bigl(r^2 \!-\! \trace{\Sigma_0} \bigr) \!+\! \trace{X} |\, \eqref{eq: LMI 1}, \eqref{eq: LMI 2}  \bigr\} .
    \end{equation}
\end{lem}
\vspace{-8pt}
\begin{pf}
    The congruence transformation of the block symmetric matrix in \eqref{eq: LMI 3} via $\Tcal$ yields that \eqref{eq: LMI 3} is equivalent to the linear matrix inequality (LMI)
    \begin{equation} \label{eq: congruence transformation for LMI 3}
        \hspace{-10pt}\begin{pmatrix}
            \beta\!-\!\gamma \|\mu_0\|^2 & * & * \\
            0 & \gamma I & * \\
            \!-g\!-\!(KCG\!-\!K^\star  \,~ K)\mu_0 \!& (KCG\!-\!K^\star  \,~ K) &\! D^{-1}\!\!
        \end{pmatrix}\!\succeq\! 0, 
    \end{equation}
where $ \Tcal:=\bigg(\!\!\!    
         \begin{smallmatrix}
            1 & 0 & 0\\
            -\mu_0 & I\, & 0\\
            0 & 0 & I
         \end{smallmatrix}\!\bigg)$ which is square and nonsingular.
It is clear that, under constraint \eqref{eq: LMI 1}, the condition that there exists a vector $g\in\Rset^{N_u}$ such that the LMI \eqref{eq: congruence transformation for LMI 3} holds is equivalent to the condition $\beta-\gamma \|\mu_0\|^2 \geq 0$. Also, a particular solution for variable $g$ to \eqref{eq: congruence transformation for LMI 3} is given by
\[g=-(KCG\!-\!K^\star \, ~ K)\mu_0.\]
Then, performing the minimisation over $\beta$ subject to the constraint $\beta-\gamma \|\mu_0\|^2 \geq 0$ proves the desired result.
\end{pf}

\begin{rem}\label{remark: reconstruct g and beta}
    Lemma \ref{lemma: eliminate g and beta} and its proof also imply that, given an optimal (feasible) solution, denoted by $(K^o, \gamma^o, X^o)$, to problem \eqref{problem: exact sdp reformulation of drro control: without g and beta}, $(K^o, -(K^oCG\!-\!K^\star \, ~ K^o)\mu_0, \gamma^o, \gamma^o \|\mu_0\|^2, X^o)$ is an optimal (feasible) solution to problem \eqref{eq:exact sdp reformulation of drro control: original}. Thus, it is preferred to solve \eqref{problem: exact sdp reformulation of drro control: without g and beta} over \eqref{eq:exact sdp reformulation of drro control: original} to obtain a DRRO controller, which is a linear controller if $\dist_0$ has a zero mean.
\end{rem}

\subsection{Elimination of the Decision Variable $K$}\label{sec: complexity reduction}
We next state a number of results as building blocks that enable the elimination of variable $K$ from problem \eqref{problem: exact sdp reformulation of drro control: without g and beta}.
\begin{rem}\label{remark:strict LMI}
    The optimal value of problem \eqref{problem: exact sdp reformulation of drro control: without g and beta} remains unchanged if the LMI \eqref{eq: LMI 2} is enforced with strict inequality.
\end{rem}

Then, by enforcing the LMI \eqref{eq: LMI 2} with strict inequality, we can apply the Projection Lemma (see, e.g. \citep[Lemma 3.1]{gahinet1994linear}, \citep[\S 2.6.2]{boyd1994linear}, and \citep[Lemma 7]{scherer1995Trends_in_Control}), which is presented as follows, and its generalisation.

\begin{lem}\label{lemma:unstructured projection lemma}
    Given matrices $U\in\Rset^{m\times n}$, $V\in\Rset^{k \times n}$, and $P\in\Sset^n$, there exists a matrix $Y\in\Rset^{m\times k}$ satisfying  
    \begin{equation}
        U^\top Y V + V^\top Y^\top U + P \succ 0, \label{eq: a general lmi in Y}
    \end{equation}
    if and only if it holds that $z^\top P z >0$ for all vectors $z \in (\ker(U) \cup \ker(V)) \!\setminus\! \{0\}$. Moreover, this condition is equivalent to that 
    \begin{equation*}
        U_\perp^\top P U_\perp \succ 0 ~\text{ and }~  V_\perp^\top P V_\perp \succ 0, 
    \end{equation*}
    where $U_\perp$ and $V_\perp$ are matrices whose columns form bases for $\ker(U)$ and $\ker(V)$, respectively.
\end{lem}

Reconstructing an unstructured solution to \eqref{eq: a general lmi in Y} involves, for example, basis extension and block matrix inversion, and efficient algorithms are available. We refer the reader to \citep[\S 7]{gahinet1994linear} for a detailed discussion. Yet, this result is not directly applicable as the variable $K$ is structured in nature.   

Building on the Projection Lemma, \citet[Theorem 2]{SCHERER1995327} provides block triangular solutions to \eqref{eq: a general lmi in Y}. We use a simple example to show its main idea as follows.

\begin{lem}\label{lemma:triangular projection lemma}
 There exists a solution $Y=\begin{pmatrix}
     Y_{11} & 0\\
     Y_{21} & Y_{22}
 \end{pmatrix}$ to \eqref{eq: a general lmi in Y}, if and only if 
 \begin{equation}\label{eq:solvability test in triangular projection lemma}
 \hspace{-1mm}    U_\perp^\top P U_\perp \!\succ 0,~ V_\perp^\top P V_\perp \!\succ 0,~\text{and}~   
     \begin{pmatrix}
         V_1 \\ U_2
     \end{pmatrix}_\perp^\top \!P \begin{pmatrix}
         V_1 \\ U_2
     \end{pmatrix}_\perp \!\!\succ 0, 
 \end{equation}
where $(U_1^\top \ \, U_2^\top)^\top =U$ and $(V_1^\top \ \, V_2^\top)^\top =V$.
\end{lem}
\vspace{-8pt}
\begin{pf}
    Given the structured variable, we rewrite \eqref{eq: a general lmi in Y} as 
    \begin{equation}\label{ineq: rewrite LMI in Y as LMI in Y1 and Y2}
        U^\top Y_1 V_1 + V_1^\top Y_1^\top U + \bar{P} \succ 0,
    \end{equation}
    where $\bar{P}:=U_2^\top Y_2 V_2 + V_2^\top Y_2^\top U_2 + P$, $Y_1:=(Y_{11}^\top \ \, Y_{21}^\top)^\top$, and $Y_2:=Y_{22}$. Then, the necessity of \eqref{eq:solvability test in triangular projection lemma} is clear. To show its sufficiency, we re-express the inequalities in \eqref{eq:solvability test in triangular projection lemma} equivalently as: $(a)~ z^\top P z>0$, $\forall z\in\ker(U)\!\setminus\!\{0\}$; $(b)~z^\top P z>0$, $\forall z\in \ker(V)\!\setminus\!\{0\}$; and $(c)~z^\top P z>0$, $\forall z \in (\ker(V_1)\cap\ker(U_2))\!\setminus\!\{0\}$, respectively. 
We first claim there exists a matrix $Y_2$ satisfying 
\begin{multline}
     ( U_2 (V_1)_\perp )^\top Y_2 ( V_2(V_1)_\perp ) + ( V_2(V_1)_\perp )^\top Y_2^\top ( U_2 (V_1)_\perp ) \\ + (V_1)_\perp^\top P (V_1)_\perp \succ 0. \label{ineq: LMI in Y_2}
\end{multline}
For any $z\neq0$ such that $V_2(V_1)_\perp z=0$, it holds that $(V_1)_\perp z \neq 0$ and $(V_1)_\perp z \in \ker(V_2)\cap \ker(V_1)$, and, from condition $(b)$, it follows that $ z^\top (V_1)_\perp^\top P (V_1)_\perp z>0 $. Similarly, for any $z\neq 0$ such that $U_2 (V_1)_\perp z=0$, it holds that $(V_1)_\perp z\neq 0$ and $(V_1)_\perp z \in \ker(U_2)\cap \ker(V_1)$, and, from condition $(c)$, it follows that $ z^\top (V_1)_\perp^\top  P (V_1)_\perp z >0$. Thus, by Lemma \ref{lemma:unstructured projection lemma}, the aforementioned claim is valid, and we choose such a matrix $Y_2$ that \eqref{ineq: LMI in Y_2} holds. We next claim there exists a matrix $Y_1$ satisfying \eqref{ineq: rewrite LMI in Y as LMI in Y1 and Y2}. For any $z\neq 0$ such that $Uz=0$, it holds that $z^\top \bar{P} z = z^\top P z>0$ by condition $(a)$. Moreover, \eqref{ineq: LMI in Y_2} already establishes that $z^\top \bar{P} z>0$ for all $z\in\ker(V_1)\!\setminus\!\{0\}$. Therefore, by Lemma \ref{lemma:unstructured projection lemma}, this proves the second claim and completes the proof.
\end{pf}

As demonstrated by the proof of Lemma \ref{lemma:triangular projection lemma}, the reconstruction of block triangular solutions needs to be implemented iteratively backwards, i.e. for $Y_2$ and then for $Y_1$. At each iteration, an LMI in the form of \eqref{eq: a general lmi in Y} is solved to obtain an unstructured solution. 

We are now prepared to present the main theorem.

\begin{thm}\label{thm: main result relying on triangular projection lemma}
Let $\bar{R}_i:=R_i (
        0_{N_y\times (N_x+N_y)} \ \, CG \ \, I \ \, 0_{N_y\times N_u} )$ and $\bar{L}_i:=L_i^\top (0_{N_u\times 2(N_x+N_y)} \  I)$ for $i=1,\ldots,T$; $\bar{R}_0:=0$ and $\bar{L}_{T+1}:=0$.
The optimal value of problem \eqref{eq:exact sdp reformulation of drro control: original} is equal to the optimal value of the following SDP
        \begin{align}
        &\inf_{ \gamma\geq0,\, X\in\Sset^{N_x\!+\!N_y} }  \gamma \bigl(r^2 \!-\! \trace{\Sigma_0} \bigr) \!+\! \trace{X}  \label{eq: final SDP}\\
        &~~~~~~~\, \text{s.t. } B_i^\top Q(\gamma,X,0) \, B_i \succ 0, \quad \forall i=1,\ldots, T+1, \nonumber
    \end{align}
    where $B_i$ is a basis matrix of the subspace $\bigcap_{j=0}^{i-1}\ker(\bar{R}_j)\cap\ker(\bar{L}_i)$.
\end{thm}
\vspace{-8pt}
\begin{pf}
By Lemma \ref{lemma: eliminate g and beta} and Remark \ref{remark:strict LMI}, the optimal value of problem \eqref{eq:exact sdp reformulation of drro control: original} is equal to 
\[\inf_{ \gamma\geq0, X\in\Sset^{N_x\!+\!N_y}\!, K\in\Kcal } \! \{ \gamma \bigl(r^2 \!-\! \trace{\Sigma_0} \bigr) +\! \trace{X} |\, Q(\gamma,X,K)\!\succ\! 0 \}, \]
where the LMI (i.e. $Q(\gamma,X,K)\!\succ\! 0$) can be rewritten in terms of the unstructured decision variables $K_i$ as
\begin{equation}\label{eq: structured LMI in summation form}
    \sum_{i=1}^T \left(\bar{L}_i^\top K_i \bar{R}_i + \bar{R}_i^\top K_i^\top \bar{L}_i \right) + Q(\gamma,X,0) \succ 0.
\end{equation}
Due to the block triangular structure of $K$, the matrices $\bar{L}_i$ naturally satisfy the condition that $\ker(\bar{L}_j) \subseteq \ker(\bar{L}_{j+1})$ for all $j=1,\ldots,T\!-\!1$, as required in \citep[Lemma 5]{scherer2013structured}. Then, by this result, there exist matrices $K_i$ satisfying the inequality \eqref{eq: structured LMI in summation form}, if and only if the $T\!+\!1$ LMIs in problem \eqref{eq: final SDP} hold, which shows the desired result. 
\end{pf}

Due to strict inequalities, only almost optimal solutions to problem \eqref{eq: final SDP} can be obtained, while this makes no significant difference in practice. 
Given an $\epsilon$-optimal solution $(\gamma',X')$ to \eqref{eq: final SDP}, one can obtain unstructured matrices $K_T'$, $K_{T-1}'$, $\ldots$, and $K_1'$ by solving $T$ LMIs in the form of \eqref{eq: a general lmi in Y} sequentially such that they satisfy \eqref{eq: structured LMI in summation form}. We refer the reader to \citep[Lemma 5]{scherer2013structured} for a detailed procedure describing this sequential reconstruction. Then, the corresponding block triangular matrix $K'$ necessarily ensures that $Q(\gamma',X',K')\succ 0$, and $(\gamma',X',K')$ is an $\epsilon$-optimal solution to problem \eqref{problem: exact sdp reformulation of drro control: without g and beta}.
\begin{rem}
Problem \eqref{eq: final SDP} naturally admits an equivalent distributed reformulation by introducing multiple copies of $(\gamma, X)$ and consensus conditions in order to decouple the $T\!+\!1$ LMIs. Thus, the objective function of such a distributed reformulation is the sum of local objective functions and also separable. This allows the implementation of consensus-based distributed algorithms (see, e.g. \citep{nedic2009cooperative,li2021distributed}), where sub-problems are solved independently or in parallel (if computational resources permit) and consensus mechanisms are designed to ensure convergence to a common solution.
This distributed optimisation offers more possibilities than \eqref{problem: exact sdp reformulation of drro control: without g and beta} to enhance scalability with respect to the time horizon $T$. Specifically, (i) the sub-problems in distributed optimisation have fewer decision variables and constraints than \eqref{problem: exact sdp reformulation of drro control: without g and beta}; (ii) the size of $B_i^\top Q \, B_i$ is generally smaller than that of $Q$ for all $i=1,\ldots, T\!+\!1$; (iii) the basis matrices $B_i$ are typically sparse and induce potential sparsity in $B_i^\top Q \, B_i$, which can be exploited to increase computational efficiency; moreover, (iv) matrices $B_i$ have certain structures to facilitate reformulating the LMIs $B_i^\top Q \, B_i\succ 0$ via Schur complements, thereby simplifying the sub-problems.
\end{rem}

\section{NUMERICAL EXAMPLE}\label{sec:numerical example}
We compare DRRO control designs via problem \eqref{eq:exact sdp reformulation of drro control: original} and problem \eqref{eq: final SDP} with controller reconstruction, respectively. Consider a discrete-time mass spring system subject to disturbance and measurement noise, where
{\small
\[ A=\begin{pmatrix}
    1 & 0.01\\
    -0.5 & 1
\end{pmatrix},~~
B=\begin{pmatrix}
    0\\
    0.5
\end{pmatrix},~~
H=\begin{pmatrix}
    0 & 0\\
    0 & 1
\end{pmatrix}.
\]}
We consider a time horizon of $T=5$ and choose weighting matrices $J_x=I$ and $J_u=I$ in the quadratic cost. Assuming $50$ i.i.d. random samples are available from the true distribution, a normal distribution $\Ncal(1,I)$, of $\mbox{col}(w,v)$, we choose a nominal distribution satisfying Assumption \ref{ass: absolute continuity} with empirical mean and covariance computed from these data and let the radius $r=\sqrt{2T}$. Problems \eqref{eq:exact sdp reformulation of drro control: original} and \eqref{eq: final SDP} are solved in \texttt{YALMIP} \citep{Lofberg2004} with \texttt{SDPT3} \citep{tutuncu2003solving}. Based on solutions to \eqref{eq: final SDP}, the feedback gain matrix $K$ is reconstructed using MATLAB function \texttt{basiclmi} iteratively. As shown in Table \ref{tab:1}, the optimal value of problem \eqref{eq: final SDP} returned by the solver is only slightly larger than that of problem \eqref{eq:exact sdp reformulation of drro control: original}, indicating negligible performance loss of the reconstructed controller. Moreover, the total computation time (averaged over 50 runs) for designing a controller in two stages is less than half of that required for solving \eqref{eq:exact sdp reformulation of drro control: original}, where the computation time for solving optimisation is the solver time returned by \texttt{YALMIP}.
\vspace{-8pt}

\begin{table}[h]
\scriptsize
\centering
\caption{Performance Comparison}
\label{tab:1}
\vspace{-3pt} 

\begin{tabularx}{0.485\textwidth}{
  >{\centering\arraybackslash}X
  >{\centering\arraybackslash}X
  >{\centering\arraybackslash}X
}
\toprule
& Optimal Value & Avg. Comp. Time \\
\addlinespace[2pt] 
\eqref{eq:exact sdp reformulation of drro control: original} & 138.476 & 0.8982 s \\
\midrule
\eqref{eq: final SDP} with reconstruction of $(K,g)$ & 138.476+$2.7\times 10^{-6}$ & 0.4142+0.0177 s \\
\bottomrule
\end{tabularx}
\end{table}

\section{CONCLUSION}
We design output feedback controllers to minimise the worst-case expected regret over all distributions in a Wasserstein ambiguity set, while this requires solving a minimax problem. Deriving strong duality results for worst-case expectation problems with quadratic objectives, we reformulate the minimax problem as an SDP. Utilising the Projection Lemma and its triangular generalisation, we eliminate certain decision variables from this SDP. The resulting problem naturally admits an equivalent distributed reformulation, enabling the deployment of distributed algorithms and with potential to improve scalability. Future work will investigate the implementation of a consensus-based distributed algorithm for this distributed SDP. 

\bibliography{reference}             
 
\end{document}